\documentclass[11pt,a4paper]{article}

\usepackage[T1]{fontenc}
\usepackage[utf8]{inputenc}
\usepackage{lmodern}
\usepackage[margin=1in]{geometry}
\usepackage{amsmath,amssymb,amsthm}
\usepackage{nicefrac}
\usepackage{graphicx}
\usepackage[authoryear,longnamesfirst]{natbib}
\usepackage{authblk}
\usepackage[colorlinks=true,linkcolor=blue,citecolor=blue,urlcolor=blue]{hyperref}

\newcommand{\udots}{\mathinner{\reflectbox{$\ddots$}}}

\newcommand{\cz}{\mathbb{C}}
\newcommand{\rz}{\mathbb{R}}
\newcommand{\im}{\mathtt{i}}

\newtheorem{theorem}{Theorem}

\title{Reduced matrix pencils for trigonometric, hyperbolic and sinc spectral estimation\thanks{This research received funding from the European Union's Horizon 2020 Research and Innovation Programme under the Marie Sk\l odowska-Curie grant agreement 101008231 (EXPOWER).}}

\author[1,2]{Annie Cuyt}
\author[1]{Wen-shin Lee\thanks{Corresponding author. ORCID: \href{https://orcid.org/0000-0002-2808-3739}{0000-0002-2808-3739}}}
\affil[1]{School of Computing, Data and Mathematical Sciences, University of Stirling, Stirling FK9 4LA, Scotland, UK}
\affil[2]{Department of Electrical Engineering, Eindhoven University of Technology, 5612 AZ Eindhoven, Netherlands}
\affil[ ]{\texttt{annie.cuyt@cemath.net}, \texttt{wen-shin.lee@stir.ac.uk}}
\date{}

\begin{document}

\maketitle

\begin{abstract}
We present $n$-sized matrix pencil methods for the spectral analysis of 
$n$-term trigonometric, hyperbolic and sinc linear sums, a problem which
occurs frequently in vibration analysis, digital signal and image
processing, direction of arrival and radio astronomy.
The new matrix pencils replace the traditional $2n$-sized matrix 
pencils used in exponential spectral analysis, which result from
re-expressing the considered functions as sums or
differences of exponentials.

The matrices in
the linear pencils are filled using $2n$ samples of the trigonometric,
hyperbolic or sinc linear sums. The frequency content is extracted
from the samples by solving a generalised eigenvalue problem,
and the linear coefficients in the
sums are computed from the right generalised eigenvectors. The latter
option is presented for each of the trigonometric, hyperbolic and sinc
cases. We numerically illustrate the novel method with an 
example from texture analysis.
\end{abstract}



\section{Introduction}

Exponential analysis, or the fitting of equidistantly sampled
data by a function of the form 
\begin{equation}
\label{1}
\sum_{j=1}^n \alpha_j \exp(\phi_j x), \qquad \alpha_j,  \phi_j \in \cz,
\quad x \in \rz, 
\end{equation}
essentially dates back to the late 18th century,
when de Prony published his exponential interpolation method in
\citet{dePr:ess:95}. In the 20th century a matrix pencil algorithm 
was published \citep{Hu.Sa:mat:90} for the
identification of the nonlinear parameters $\phi_j$ from generalised
eigenvalues of a linear pencil of $n \times n$ matrices.
The linear parameters $\alpha_j$ were recovered in a subsequent step
by solving a Vandermonde-structured linear system. We show here how they can also be
obtained from the computation of the eigenvectors of the pencil.

In the meantime exponential analysis was
generalised to functions other than the exponential or the power function \citep{Pe.Pl:gen:13,Pl.St.ea:rec:19},
and to higher dimensions, but often without an accompanying 
matrix pencil version. The basics from 
\citet{dePr:ess:95} and \citet{Hu.Sa:mat:90} are recapped in Section 2. 
In the subsequent sections a matrix pencil algorithm is developed for
parametric models of the type
\begin{gather}
\sum_{j=1}^{n_c} \alpha_{j1} \cos(\phi_{j1} x) + \sum_{j=1}^{n_s}
\alpha_{j2} \sin(\phi_{j2} x),  \label{CoSi} \\
\sum_{j=1}^{n_c} \alpha_{j1} \cosh(\phi_{j1} x) + \sum_{j=1}^{n_s}
\alpha_{j2} \sinh(\phi_{j2} x), \label{ChSh} \\
\sum_{j=1}^n \alpha_j \text{sinc}(\phi_j x) \label{Sinc} 
\end{gather}
and their multivariate counterparts, where all linear and nonlinear
parameters are real.
How to include phase shifts in all of the above models is indicated when
discussing the trigonometric generalisation. 

Note that with $n_c, n_s$ and $n$ indicating respectively the numbers of
nonzero distinct terms in \eqref{CoSi}, \eqref{ChSh} and \eqref{Sinc},
it is not only tacitly assumed that the
$\alpha_j, \alpha_{j1}$ and $\alpha_{j2}$ are all nonzero and that the
$\phi_{j1}$ and the $\phi_{j2}$ are mutually distinct, as in the
exponential case of Section 2, but also that 
$\phi_j \not= \pm\phi_i$ for $j \not= i$
and that in the odd parts $\phi_j\not=0$.
Otherwise certain (hyperbolic) cosine or (cardinal) sine 
terms would collide or vanish and certain matrices in the subsequent theorems 
would become singular.
This assumption therefore holds throughout Sections 3, 4 and 5,
wherever required.

\section{Matrix pencil method for exponential analysis}

Let us consider the nonlinear interpolation problem
\begin{equation}
\label{sidi}
\sum_{j=1}^n \alpha_j \exp(\phi_j x_k) = f_k, \qquad k=0, \dots, 2n-1,
\end{equation}
where $x_k=k\Delta$ with $\Delta\neq0$ and the $\phi_j$ are 
mutually distinct and satisfy the Nyquist condition \citep{Ny:cer:28}
\begin{equation}
\label{nyquist}
\max_{j=1, \ldots, n}|\Im(\phi_j)\Delta| < \pi,
\end{equation}
where $\Im(\cdot)$ denotes the imaginary part of a complex number. In the
transition to real models in the subsequent sections, the Nyquist
constraint becomes $\max_{j=1, \ldots, n} |\phi_j\Delta| < \pi$.
In \eqref{sidi} the $x_k$ and $f_k$ are given and the $\alpha_j$ and
$\phi_j$ are to be determined.
Often also $n$ is unknown.

Exponential analysis is actually a reconstruction problem rather than an
approximation technique. The derivation of what is basically Prony's
method in \citet{dePr:ess:95} 
is based on the assumption that the data $f_k$ are taken from
a function $f(x)$ that is of the form \eqref{1}.
Then the unknown parameters $\alpha_j$ and $\phi_j$ 
can be extracted from such data $f_k$ as follows.
If we denote $\Phi_j := \exp(\phi_j \Delta)$,
then it is clear that the data $f_k$ are structured as
\begin{equation}
\label{S}
f_k = \sum_{j=1}^n \alpha_j \Phi_j^k, \qquad k=0,\dots,2n-1.
\end{equation}
We now want to obtain the values $\Phi_j, j=1, \ldots, n$ and $\alpha_j,
j=1, \ldots, n$ from the $2n$ samples $f_k$. From $\Phi_j$ the unknown
$\phi_j$ can easily be deduced because
$|\Im(\phi_j) \Delta|<\pi$ and hence no periodicity problem arises. 
Also, since the $\phi_j$ are mutually distinct, the $\Phi_j$ are mutually
distinct and hence the Vandermonde matrix in \eqref{F} below is regular.
In addition to computing $\Phi_j$ and $\alpha_j$, we indicate what can be
deduced about $n$. Knowledge of the latter is crucial for the success of the
algorithm. Temporarily we assume that $n$ is known.

Let us denote by $H_n^{(r)}$ the Hankel matrix
\begin{equation} \label{Hrevision}
H_n^{(r)} := 
\begin{pmatrix}
f_r & f_{r+1} &\dots & f_{r+n-1} \\
f_{r+1} & & & \\
\vdots & \udots & & \vdots \\
f_{r+n-1} & & \dots & f_{r+2n-2}
\end{pmatrix}.
\end{equation}
From the matrix factorisations
\begin{equation}
\begin{aligned}
\label{F}
H_n^{(0)} &= V_n D_\alpha V_n^T,  \\ 
H_n^{(1)} &= V_n D_\alpha D_\Phi V_n^T, 
\end{aligned}
\end{equation}
where $V_n, D_\alpha$ and $D_\Phi$ respectively denote the square
Vandermonde matrix
\begin{equation*}
V_n = \begin{pmatrix}
1 & 1 & \dots & 1 \\
\Phi_1 & \Phi_2 & \dots & \Phi_n \\
\vdots & \vdots & & \vdots \\
\Phi_1^{n-1} & \Phi_2^{n-1} & \dots & \Phi_n^{n-1}
\end{pmatrix}
\end{equation*}
and the diagonal matrices
$D_\alpha=\text{diag} (\alpha_1, \ldots, \alpha_n)$, 
$D_\Phi = \text{diag}(\Phi_1, \ldots, \Phi_n)$,
it is straightforward to see that the $\Phi_j$ can be 
obtained as generalised eigenvalues \citep{Hu.Sa:mat:90,Go.Mi.ea:sta:99}. 
The $\Phi_j$ satisfy
\begin{equation}
\label{G}
\det \left( H_n^{(1)} - \Phi_j H_n^{(0)} \right) = 0, \qquad j=1,
\ldots, n. 
\end{equation}
The coefficients $\alpha_j$ in the model \eqref{1} can be obtained from
any set of $n$ interpolation conditions taken from \eqref{S}. 
With $\Phi_1, \ldots, \Phi_n$ computed from \eqref{G}, 
the remaining equations in \eqref{S} are linearly dependent. However,
usually the full system \eqref{S} is solved in the least-squares sense.
Another way to obtain the linear coefficients $\alpha_j$ is the following. 
Let $v_j$ be a right eigenvector of the generalised 
eigenvalue problem \eqref{G} and let $I_n$ be the $n \times n$ identity matrix
and $e_j$ the vector with 1 at position $j$ and zero elsewhere.
Then from \eqref{F} and \eqref{G} we see that 
\begin{equation}
V_n D_\alpha \left( D_\Phi - \Phi_j I_n \right) V_n^T v_j = 0, \qquad j=1,
\ldots, n
\end{equation}
and thus that, with the matrix $D_\Phi - \Phi_j I_n$ being diagonal, 
the vector $V_n^T v_j$ is a multiple of $e_j$ or
$$v_j = \gamma_j V_n^{-T}e_j, \qquad j=1, \ldots, n,$$
where $V_n^T (v_1, \ldots, v_n) = \text{diag}(\gamma_1, \ldots, \gamma_n)$. 
Consequently
$$H_n^{(0)}v_j=\alpha_j \gamma_j V_n e_j, \qquad j=1, \ldots, n$$
from which we can obtain and average
$n$ estimates for the $\alpha_j$ as
$$\alpha_j = \text{mean}_{i=1, \ldots, n}\; {\left( H_n^{(0)}v_j \right)_i
\over \gamma_j \Phi_j^{i-1}}, \qquad j=1, \ldots, n.$$

Now what can be said about $n$? Let $|\cdot|$ denote the determinant 
det$(\cdot)$ of a matrix. In the noise-free case, its value can
be nailed down quite precisely. We
read in \citet{He:app:74} and \citet{Ba.Gr:pad:96} that on the one hand,
for $N<n$ and $r\ge
0$, $|H^{(r)}_N|$ is only accidentally zero, depending on the value of
$\Delta$, while on the other hand, for $N>n$ and $r \ge 0$, $|H^{(r)}_N|$
is always zero, irrespective of the value of $\Delta$. 
Most importantly, for $N=n, r \ge 0$ and $\Phi_i \ne \Phi_j$
when $i \ne j$, one has $|H^{(r)}_n| \neq 0$ whatever the value of $\Delta$.
In order to inspect $|H_N^{(r)}|$ for $N>n$, additional samples up to
$f_{r+2N-2}$ need to be provided, in other words at least the additional
sample $f_{2n}$ (in case $r=0$ and $N=n+1$). 
A nice discussion, based on algebraic arguments, is presented in 
\citet{Ka.Le:ear:03}. How $n$ can reliably be determined when working with
noisy data is discussed in \citet{Br.Cu.ea:vex:20}.

\section{Matrix pencil method for the trigonometric model}

We now consider the reconstruction of functions of
the form \eqref{CoSi}. 
If either $n_s=0$ or $n_c=0$, then it is sufficient to collect data at
respectively $x_k=k\Delta$ for $k=0, \ldots, 2n-1$ or at 
$x_k=k\Delta$ for $k=1, \ldots, 2n$, where $n=\max(n_c, n_s)$, because
$\cos(-\phi_j k\Delta) = \cos(\phi_j k\Delta)$, $\sin(-\phi_j
k\Delta) = -\sin(\phi_j k\Delta)$ and $\sin(0)=0$. Since we assumed that
$\phi_j \not= \pm\phi_i$ for $j \not= i$ and that in
the odd part $\phi_j \not=0$,
we know that under the Nyquist constraint
\eqref{nyquist}, the absolute values $0 \le |\phi_j\Delta|<\pi$ are
mutually distinct and that $\phi_j\Delta=0$ is excluded in sine sums.

When $n_c n_s \not= 0$ then $n=n_c$ if $n_c>n_s$ and $n=n_s$ if
$n_c \le n_s$. In the former case, we
collect the data $f_k$ by sampling 
\eqref{CoSi} at $x_k=k\Delta$ with $k=-2n+1, -2n+2, \ldots, 0, \ldots,
2n-2, 2n-1$ and in the latter case by sampling at $x_k=k\Delta$ for
$k=-2n, \ldots, -1, 0, 1, \ldots, 2n$. 
With these samples we then separate the even and odd parts of
\eqref{CoSi} by rewriting
\begin{align*}
f_k &= {f_k+f_{-k} \over 2} + {f_k-f_{-k} \over 2}, \qquad
k=0, \ldots, 2n-1, \\
&= \sum_{j=1}^{n_c} \alpha_{j1}\cos(\phi_{j1}k\Delta) + 
\sum_{j=1}^{n_s} \alpha_{j2} \sin(\phi_{j2}k\Delta).
\end{align*}
Of course, in general, both $n_c$ and $n_s$ are unknown and so, in practice,
$n$ is overestimated. Then the first sampling scheme can be used in all
situations. After obtaining the samples $(f_k+f_{-k})/2$ of the even part
and the samples $(f_k-f_{-k})/2$ of the odd part, we can 
proceed with the  of \eqref{CoSi} where respectively
$n_s=0$ and $n_c=0$. 

In both cases, from here on, 
we call the unknown coefficients
$\alpha_j$, the unknown frequencies $\phi_j$ and the samples of the even
or odd part $g_k$ with the number of terms being denoted by $n$.
We reserve the notation $f_k$ for samples directly taken from
\eqref{CoSi}, \eqref{ChSh} or \eqref{Sinc}.
Before doing so, 
we point out how the new method differs from the mere interpolation
of \eqref{CoSi} by $2n_c+2n_s$ complex conjugate exponentials. 

First of all, instead of requiring $4n_c+4n_s$ samples to fit the
complex conjugate exponential terms, only half as many samples are
needed, more precisely $4\max(n_c, n_s) \mp 1$.
Also, the matrix size of the structured generalised eigenvalue problems 
is reduced 
from $2(n_c + n_s) \times 2(n_c+n_s)$ for the exponential model with
complex conjugate terms,
to $n_c \times n_c$ for the even part and $n_s \times n_s$ for the odd part,
which reduces the computational complexity.

\subsection{The cosine case: even part ($n=n_c$)}

Besides the Hankel-structured matrix $H_n^{(r)}$ in \eqref{Hrevision}, 
now filled with samples $g_k$
from the even part of \eqref{CoSi}, we introduce the 
following Hankel and Toeplitz-structured matrices:
\begin{align*}
H_{n,+}^{(r)} &:= H_n^{(r)}, \\
H_{n,-}^{(r)} &:= \begin{pmatrix}
g_r & g_{r-1} & \dots & g_{r-n+1} \\
g_{r-1} & & & \\
\vdots & \udots & & \vdots \\
g_{r-n+1} & & \dots & g_{r-2n+2}
\end{pmatrix}
\end{align*}
and 
\begin{align*}
T_n^{(r)} &:= \begin{pmatrix} g_r & g_{r-1} &
\cdots & g_{r-n+1} \\
g_{r+1} & & & \\
\vdots & & \ddots & \vdots \\ 
g_{r+n-1} & & \cdots & g_r
\end{pmatrix} , \\
T_{n,+}^{(r)} &:= T_n^{(r)}, \qquad
T_{n,-}^{(r)} := \left( T_n^{(r)} \right)^T.
\end{align*}
For a pure cosine model and with $r=0$ the Toeplitz-structured matrices
are symmetric.
Now consider the structured matrix
\begin{equation*}
C_n^{(r)} := {1 \over 4} H_{n,+}^{(r)} + 
{1 \over 4} H_{n,-}^{(r)} + 
{1 \over 4} T_{n,+}^{(r)} +
{1 \over 4} T_{n,-}^{(r)}.
\end{equation*}
When $r=0$, the first two matrices in the sum
coincide and the latter two do as well.

We emphasise that working directly with the cosine sum, 
instead of expressing
it in terms of exponentials as $\cos x = (\exp(\im x) + \exp(-\im
x))/2$, allows us to reduce the size of the matrices involved 
in the generalised eigenvalue problems.

\begin{theorem}
The structured matrix $C_n^{(r)}$ factorises as
\begin{align*}
C_n^{(r)} &= W_n L^{(r)}_n D_\alpha W_n^T, \\ 
W_n &=\begin{pmatrix}
1 & \cdots  & 1 \\
\cos(\phi_1\Delta) & \cdots & \cos(\phi_n\Delta) \\
\vdots & & \vdots \\
\cos(\phi_1(n-1)\Delta) & \cdots & \cos(\phi_n(n-1)\Delta) 
\end{pmatrix}, \\ 
L^{(r)}_n &= \text{\rm diag}(\cos(\phi_1 r\Delta), \ldots, \cos(\phi_n r\Delta)).
\end{align*}
\end{theorem}

\begin{proof}
The proof is a mere verification of the matrix product entry at position
$(k, \ell)$ for $k, \ell=1, \ldots, n$ in combination with the use of the
trigonometric sum formulas. Note that the matrix $W_n$ 
is regular because the $|\phi_j\Delta|$ are
mutually distinct.
We have 
$$|W_n| = 2^{(n-1)(n-2)/2} \prod_{1\le i<j\le n} \left(\cos(\phi_j\Delta)
-\cos(\phi_i\Delta) \right).$$
\end{proof}

From this matrix factorisation we deduce that for
$\lambda_j=\cos(\phi_j\Delta)$,
\begin{equation*}
\det \left( C_n^{(1)} - \lambda_j C_n^{(0)} \right) = \det \left(
W_n \left( L_n^{(1)} - \lambda_j L_n^{(0)} \right) D_\alpha W_n^T \right)=0,
\end{equation*}
which opens the
door to the generalised eigenvalue problem formulation
\begin{equation}
\label{gep_c}
\left( C_n^{(1)} \right) v_j = \cos(\phi_j \Delta) 
\left( C_n^{(0)} \right) v_j, \qquad j=1, \ldots, n, 
\end{equation}
where the $v_j$ are generalised right eigenvectors. Because $g_k=g_{-k}$,
setting up \eqref{gep_c} only takes the $2n$ samples $g_k, k=0, \ldots, 2n-1$. 
Under condition \eqref{nyquist}, 
the $\phi_j$ can unambiguously be extracted from the generalised eigenvalues
$\cos(\phi_j\Delta)$, up to an irrelevant sign. 

The coefficients $\alpha_j, j=1, \ldots, n$ are computed from an $n \times n$
linear subsystem of the interpolation conditions
\begin{equation*}
\sum_{j=1}^n \alpha_j \cos(\phi_j k\Delta) = f_k, \qquad 
k=0, \ldots, 2n-1.
\end{equation*}
To extract the $\alpha_j$ from right eigenvectors $v_j$, we first obtain, as
in Section 2, that $W_n^T v_j$ is a multiple of $e_j$ and thus conclude
that $C^{(0)}_n v_j$ equals
$$C^{(0)}_n v_j = \alpha_j \gamma_j W_n e_j, \qquad j=1, \ldots, n$$
with $W_n^T(v_1, \ldots, v_n) = \text{diag}(\gamma_1, \ldots, \gamma_n)$.
This gives the formula
$$\alpha_j = \text{mean}_{i=1, \ldots, n}\; {\left( C_n^{(0)}v_j \right)_i
\over \gamma_j \cos((i-1)\phi_j\Delta)}.$$

\subsection{The sine case: odd part ($n=n_s$)}

We now define
\begin{equation*}
S_n^{(r)} := {1 \over 4} H_n^{(1+r)} + 
{1 \over 4} H_n^{(1-r)} + 
{1 \over 4} T_n^{(1+r)} + 
{1 \over 4} T_n^{(1-r)},
\end{equation*}
where all matrices are filled with the samples $g_k$ from the odd part of
\eqref{CoSi}.

\begin{theorem}
The structured matrix $S_n^{(r)}$ factorises as
\begin{align*}
S_n^{(r)} &= U_n L^{(r)}_n D_\alpha W_n^T, \\
U_n &= \begin{pmatrix} 
\sin(\phi_1\Delta) & \cdots & \sin(\phi_n\Delta) \\
\vdots & & \vdots \\
\sin(\phi_1 n\Delta) & \cdots & \sin(\phi_n n\Delta) 
\end{pmatrix}, \\
W_n &=\begin{pmatrix}
1 & \cdots  & 1 \\
\cos(\phi_1\Delta) & \cdots & \cos(\phi_n\Delta) \\
\vdots & & \vdots \\
\cos(\phi_1(n-1)\Delta) & \cdots & \cos(\phi_n(n-1)\Delta)
\end{pmatrix}, \\ 
L^{(r)}_n &= \text{\rm diag}(\cos(\phi_1 r\Delta), \ldots,
\cos(\phi_n r\Delta)).
\end{align*} 
\end{theorem}

\begin{proof}
The proof is again a verification of the matrix product entry, at the
position $(k, \ell)$ with $k, \ell =1, \ldots, n$, in combination with the 
use of the trigonometric sum formulas. Note that the matrices
$U_n$  and $W_n$ are regular because
the $|\phi_j\Delta|$ are mutually distinct and $\phi_j\Delta \not=0$.
We have
$$|U_n| = 2^{n(n-1)/2} \prod_{j=1}^n \sin(\phi_j\Delta) 
\prod_{1 \le i<j \le n} \left( \cos(\phi_j\Delta) - \cos(\phi_i\Delta) 
\right).$$
\end{proof}

From this decomposition we find that for $\lambda_j = \cos(\phi_j\Delta)$,
\begin{equation*}
\det \left( S_n^{(1)} - \lambda_j S_n^{(0)} \right) = \det \left( U_n
\left( L_n^{(1)} - \lambda_j L_n^{(0)} \right) D_\alpha W_n^T \right) = 0
\end{equation*}
and thus that the values $\cos(\phi_j\Delta), j=1,
\ldots, n$ are obtained as the generalised eigenvalues of the problem
\begin{equation}
\label{gep_s}
\left( S_n^{(1)} \right) v_j = \cos(\phi_j\Delta)
\left( S_n^{(0)} \right) v_j, \qquad j=1, \ldots, n. 
\end{equation}
We point out that setting up this generalised eigenvalue problem requires
samples at $x_{-n+1}, \ldots$, $x_{2n}$. Since
$g_k=-g_{-k}$ and $g_0=0$ it only costs $2n$ samples.

From the values $\cos(\phi_j\Delta)$ we can extract $\pm|\phi_j|\Delta$,
without specification of the sign though. This, however, does not cause a problem,
as the odd part of \eqref{CoSi} remains identical whether the terms in the
model are $\alpha_j \sin(\phi_j x)$ or $(-\alpha_j) \sin(-\phi_j x)$. So
after extracting the $|\phi_j|$ from  the generalised eigenvalues
$\cos(\phi_j\Delta)$, we compute the values $\pm\sin(\phi_j k\Delta)$ and
solve for the $\pm\alpha_j, j=1, \ldots, n$
from an $n \times n$ linear subsystem of the interpolation conditions
\begin{equation*}
\sum_{j=1}^n \alpha_j \sin(\phi_j k\Delta) = f_k, 
\qquad k=1, \ldots, 2n.
\end{equation*}

Generalised eigenvectors $v_j$ are again multiples of $W_n^{-T}e_j$ and so we
can also compute the linear coefficients $\alpha_j$ from
$$S_n^{(0)} v_j = \alpha_j \gamma_j U_n e_j, \qquad j=1, \ldots, n$$
where $W_n^T (v_1, \ldots, v_n) = \text{diag} (\gamma_1, \ldots,
\gamma_n)$, as
$$\alpha_j = \text{mean}_{i=1, \ldots, n}\; {\left( S_n^{(0)}v_j \right)_i
\over \gamma_j \sin(i\phi_j\Delta)}.$$

\subsection{Dealing with phase shifts}

In each of \eqref{CoSi}, \eqref{ChSh} and \eqref{Sinc} phase shifts can be
included. Without loss of generality, we here describe how to deal with 
them in the multi-sinusoidal model
\begin{equation}
\label{phase}
\sum_{j=1}^n \alpha_j \sin(\phi_j x + \omega_j), \qquad |\omega_j|<\pi/2. 
\end{equation}

When using the classical complex exponential approach, each term in
\eqref{phase} is rewritten as 
$$\frac{\alpha_j \exp(\im \omega_j) \exp(\im \phi_j x) - \alpha_j \exp(-\im
\omega_j) \exp(-\im \phi_j x)}{2 \im},$$
leading to a $(2n)$-term model, requiring $4n$ samples for the
identification of the $\alpha_j, \phi_j, \omega_j, j=1, \ldots, n$. With
the current approach, each term in \eqref{phase} is rewritten as
$$\alpha_j \cos(\omega_j) \sin(\phi_j x) + \alpha_j \sin(\omega_j)
\cos(\phi_j x).$$
Its even and odd parts respectively equal 
$\alpha_j \sin(\omega_j) \cos(\phi_j x)$ and 
$\alpha_j \cos(\omega_j) \sin(\phi_j x)$ 
and each leads to an $n$-term trigonometric problem from which the
frequencies $\phi_j$ and the linear coefficients 
$\alpha_j\cos(\omega_j)$, $\alpha_j\sin(\omega_j)$
can be identified, given either $2n$ samples of \eqref{phase}
for $x>0$ or $2n$ samples of \eqref{phase} 
for $x<0$. Then $\omega_j$ is computed from
$\tan(\omega_j)$ and $\alpha_j$ from the two extracted linear coefficients.

\section{Matrix pencil method for the hyperbolic and sinc models}

\subsection{The hyperbolic case}

For the reconstruction, through interpolation, of expression \eqref{ChSh},
one proceeds in a similar way. 
Let $n=\max(n_c, n_s)$ and sample
\eqref{ChSh} at the points $x_k=k\Delta$ for $k=-2n, \ldots, -1, 0, 1,
\ldots, 2n$. Since we are not dealing with periodic functions now,
condition \eqref{nyquist} need not be imposed, but for regularity of the
structured matrices, we maintain the assumption 
that in $W_n$ the absolute
values $|\phi_j\Delta|$ are mutually distinct and that in $U_n$ also
$\phi_j\Delta=0$ is excluded. Just as earlier, we
separate the even and odd parts of the model by constructing data
$(f_k+f_{-k})/2$ for the hyperbolic cosine part and $(f_k-f_{-k})/2$ for
the hyperbolic sine part.

With these data, the matrices $C_n^{(r)}$ and $S_n^{(r)}$, defined in the
previous section and filled with the samples $g_k$ of respectively 
the even or odd part of \eqref{ChSh}, now factorise as follows. 

\begin{theorem}
The structured matrix $C_n^{(r)}$ factorises as
\begin{align*}
C_n^{(r)} &= W_n L^{(r)}_n D_\alpha W_n^T, \\ 
W_n &=\begin{pmatrix}
1 & \cdots  & 1 \\
\cosh(\phi_1\Delta) & \cdots & \cosh(\phi_n\Delta) \\
\vdots & & \vdots \\
\cosh(\phi_1(n-1)\Delta) & \cdots & \cosh(\phi_n(n-1)\Delta) 
\end{pmatrix}, \\ 
L^{(r)}_n &= \text{\rm diag}(\cosh(\phi_1 r\Delta), \ldots, \cosh(\phi_n r\Delta)).
\end{align*}
\end{theorem}

\begin{theorem}
The structured matrix $S_n^{(r)}$ factorises as
\begin{align*}
S_n^{(r)} &= U_n L^{(r)}_n D_\alpha W_n^T, \\
U_n &= \begin{pmatrix} 
\sinh(\phi_1\Delta) & \cdots & \sinh(\phi_n\Delta) \\
\vdots & & \vdots \\
\sinh(\phi_1 n\Delta) & \cdots & \sinh(\phi_n n\Delta) 
\end{pmatrix}, \\
W_n &=\begin{pmatrix}
1 & \cdots  & 1 \\
\cosh(\phi_1\Delta) & \cdots & \cosh(\phi_n\Delta) \\
\vdots & & \vdots \\
\cosh(\phi_1(n-1)\Delta) & \cdots & \cosh(\phi_n(n-1)\Delta)
\end{pmatrix}, \\
L^{(r)}_n &= \text{\rm diag}(\cosh(\phi_1 r\Delta), \ldots,
\cosh(\phi_n r\Delta)).
\end{align*} 
\end{theorem}

\begin{proof}
In the proofs one merely needs to take into account the different
addition and subtraction rules
\begin{align*}
\cosh(A+B) &= \cosh(A) \cosh(B) + \sinh(A)\sinh(B) \\
\cosh(A-B) &= \cosh(A) \cosh(B) - \sinh(A)\sinh(B) 
\end{align*}
and
\begin{align*}
\sinh(A+B) &= \sinh(A) \cosh(B) + \cosh(A)\sinh(B) \\
\sinh(A-B) &= \sinh(A) \cosh(B) - \cosh(A)\sinh(B) . 
\end{align*}
\end{proof}

\subsection{The sinc function case}

The sinc function is classically defined as $\text{sinc}(x) =
\sin(x)/x$. Our interest is in the reconstruction 
of model \eqref{Sinc} under condition \eqref{nyquist}. 
With the samples $f_k$ of \eqref{Sinc} 
taken at $x_k=k\Delta$ for $k=1, \ldots, 2n$ we
construct the values
\begin{equation*}
F_k := k\Delta f_k, \qquad k=1, \ldots, 2n.
\end{equation*}
Then the interpolation problem
\begin{equation*}
\sum_{j=1}^n \alpha_j \text{sinc}(\phi_j x_k) = f_k, \qquad k=1, \ldots, 2n
\end{equation*}
can be rewritten as
\begin{equation}
\label{newSinc}
\sum_{j=1}^n \left( {\alpha_j \over \phi_j} \right) 
\sin (\phi_j k \Delta) = F_k, \qquad k=1, \ldots, 2n. 
\end{equation}
From this rewrite, the parameters $\phi_j$ and $\alpha_j$ can be obtained
by applying the method from Section 3.2. Since the $\phi_j$ are obtained
first and separately from the $\alpha_j$, it is straightforward to
extract the $\alpha_j$ from the coefficients $\alpha_j/\phi_j$ in
\eqref{newSinc}.

\section{Multivariate versions}

Each of the models \eqref{CoSi}, \eqref{ChSh} and \eqref{Sinc} can be 
generalised to the multivariate case, similar to 
the multivariate generalisation in \citet{Cu.Le:mul:18}
of the exponential model \eqref{sidi},
although some extra details need to be taken care of
when dealing with the base terms in \eqref{CoSi}, \eqref{ChSh} or
\eqref{Sinc}. 
Without loss of generality we present the bivariate case and explain 
the interpolation of the even part of the multivariate version of
\eqref{CoSi},
$$g(x,y) = \sum_{j=1}^n \alpha_j \cos(\phi_j x + \psi_j y),$$
leaving a similar treatment of the odd part of \eqref{CoSi}, as well as 
\eqref{ChSh} and \eqref{Sinc} to the reader. Just as in the univariate
case, the linear coefficients $\alpha_j$, the frequency components
$\phi_j$ and $\psi_j$ are unknown. The purpose is to extract them from a
minimal number of samples, in terms of a generalised eigenvalue problem.
For simplicity, we assume in what follows that $n$ is known.

Before we continue, we point
out that the analysis of the even part, when expressed in terms of $N=2n$ 
complex exponentials, requires $3N=6n$ samples. The new approach brings 
this down to only $4n$.

We define two linearly independent sample
vectors $\Delta=(\Delta_1, \Delta_2)$ and $\delta=(\delta_1, \delta_2)$
satisfying

\begin{equation}
\label{nyquist2}
\begin{gathered}
\max_{j=1, \ldots, n} |\phi_j \Delta_1 + \psi_j \Delta_2| < \pi, \\
\max_{j=1, \ldots, n} |\phi_j \delta_1 + \psi_j \delta_2| < \pi.
\end{gathered}
\end{equation}
In the multivariate generalisation, the role of the argument $\phi_j\Delta$ 
in the previous sections
is taken by the inner product 
$$P_j:= \phi_j \Delta_1 + \psi_j \Delta_2.$$
Analogously we define the inner product 
$$Q_j:= \phi_j\delta_1 + \psi_j\delta_2.$$
In the multivariate setting we assume that the $(\phi_j,\psi_j), 
j=1, \ldots, n$ are mutually distinct pairs and that no addition or
subtraction of frequency pairs gives zero.
In addition we need $\Delta=(\Delta_1, \Delta_2)$ to guarantee to be such
that the $|P_j|$ are mutually distinct.
Let us sample
\begin{align*}
g_k &= g(k\Delta_1, k\Delta_2), \qquad k= 0, \ldots, 2n-1, \\
G_{k,+} &= g(k\Delta_1+\delta_1, k\Delta_2+\delta_2), \qquad k=0, \ldots, n-1, \\
G_{k,-} &= g(k\Delta_1-\delta_1, k\Delta_2-\delta_2), \qquad k=0, \ldots, n-1,
\end{align*}
with $g_{-k}=g_k$ and $G_{0,+}=G_{0,-}$.
As in \eqref{gep_c} we find that the $\cos(P_j)$, $j=1, \ldots, n$ can
be obtained from  the generalised eigenvalue problem
\begin{equation}
\label{gep_c2}
\left( C_n^{(1)} \right) v_j = \cos(P_j)
\left( C_n^{(0)} \right) v_j, \qquad j=1, \ldots, n, 
\end{equation}
where the matrices $C_n^{(r)}$ are now filled with the bivariate samples $g_k$
in which the univariate $\phi_j\Delta$ is replaced by the bivariate $P_j$.
The linear coefficients 
$\alpha_j, j=1, \ldots, n$ are obtained from an $n \times n$ subsystem of
\begin{equation}
\label{alpha2}
\sum_{j=1}^n \alpha_j \cos(kP_j) = g_k, \qquad k=0, \ldots, 2n-1,
\end{equation}
or from the right generalised eigenvectors $v_j$ of \eqref{gep_c2}.

Because of \eqref{nyquist2}, $|P_j|$ can be extracted uniquely from
$\cos(P_j)$. Since
$$
G_{k,+} + G_{k,-} = 2 \sum_{j=1}^n \alpha_j \cos(Q_j) \cos(k P_j), \qquad
k=0, \ldots, n-1,$$
the unknowns $\cos(Q_j), j=1, \ldots, n$
are obtained from the Vandermonde-like linear system 
\begin{equation}
\label{VDM2}
{\cal V}_n
\begin{pmatrix} \alpha_1 \cos(Q_1) \\ \vdots \\ \alpha_n \cos(Q_n)
\end{pmatrix} = {1 \over 2} \begin{pmatrix} G_{0,+}+G_{0,-} \\ \vdots \\
G_{n-1,+}+G_{n-1,-} \end{pmatrix} ,
\end{equation} 
where
$${\cal V}_n = 
\begin{pmatrix} 1 & \cdots & 1 \\
\cos(P_1) & \cdots & \cos(P_n) \\ \vdots & & \vdots \\ \cos((n-1)P_1) &
\cdots & \cos((n-1)P_n) \end{pmatrix} .$$
We emphasise that, with the $|P_j|$ mutually distinct, ${\cal V}_n$ is
regular.
From $\cos(Q_j)$ again $|Q_j|, j=1, \ldots, n$ can be extracted. It remains to
discuss how to identify the parameters $\phi_j, \psi_j, j=1, \ldots, n$ in
$g(x,y)$ from the knowledge of $|P_j|, |Q_j|, j=1, \ldots, n$. At first
sight, there are four solutions from solving the $2 \times 2$ linear systems
$$\left\{ \begin{aligned} & \Delta_1 \phi_j + \Delta_2 \psi_j = \pm|P_j| \\
& \delta_1 \phi_j + \delta_2 \psi_j = \pm|Q_j| \end{aligned}
\right.,$$
unless some additional information on the sign of $P_j$ and $Q_j$ is
provided.
Since $g(x,y)$ consists of cosine terms, the combinations $(+,+)$
and $(-,-)$ deliver the same model, and similarly for $(+,-)$ and $(-,+)$.
Therefore only two possibilities remain. 

In order to identify the correct
solution between these two, we require only two more samples while 
reusing the coefficient matrix ${\cal V}_n$ at the same time.
When adding the samples $G_{-1,-}$ and $G_{n,+}$ and making use of the
identities
\begin{align*}
G_{k,+}&= \sum_{j=1}^n \alpha_j \left( \cos((k-1)P_j) \cos(P_j+Q_j) - \sin((k-1)P_j) \sin(P_j+Q_j) \right), \\
G_{k-2,-}&= \sum_{j=1}^n \alpha_j \left( \cos((k-1)P_j) \cos(P_j+Q_j)  +  \sin((k-1)P_j) \sin(P_j+Q_j) \right),
\end{align*}
then $\cos(P_j+Q_j)$ can be obtained from the linear system
$${\cal V}_n \begin{pmatrix} \alpha_1 \cos(P_1+Q_1) \\ \vdots \\ \alpha_n \cos(P_n+Q_n)
\end{pmatrix} = {1 \over 2} \begin{pmatrix} G_{1,+}+G_{-1,-} \\ \vdots \\
G_{n,+}+G_{n-2,-} \end{pmatrix} .$$ 
Since we can assume that not both $P_j$ and $Q_j$ equal zero, because otherwise
$\phi_j=0=\psi_j$, we check whether the solution of $(+,+)$ or $(+,-)$ matches 
$|P_j+Q_j|$, extracted from $\cos(P_j+Q_j)$. Note that it is impossible for
$P_j+Q_j$ to equal $P_j-Q_j$ and at the same time for $-P_j-Q_j$ to 
equal $P_j-Q_j$ unless both $P_j$ and $Q_j$ are zero.

When $n$ is small and sample collection is expensive, one may for
instance prefer to
only collect the sample $G_{n,+}$ and
not collect the sample $G_{-1,-}$, thereby losing the possibility to
reuse ${\cal V}_n$. In fact, because of the identity
\begin{equation}
\label{half}
\cos\left( {2k-1 \over 2}P_j \right) \cos \left( {2k+1 \over 2}P_j \right) 
=  {1 \over 2} \left( \cos(P_j) + \cos(2kP_j) \right),
\end{equation}
the values
$\cos({\nicefrac{1}{2}} P_j+Q_j)$ can be obtained from the linear system
\begin{equation*}
{\cal W}_n \begin{pmatrix} \alpha_1 \cos(\nicefrac{1}{2}P_1+Q_1) \\ \vdots
\\ \alpha_n \cos(\nicefrac{1}{2}P_n+Q_n)
\end{pmatrix} = {1 \over 2} \begin{pmatrix} G_{1,+}+G_{0,-} \\ \vdots \\
G_{n,+}+G_{n-1,-} \end{pmatrix},
\end{equation*}
where 
$${\cal W}_n = \begin{pmatrix}  \cos(\nicefrac{1}{2}P_1) & \cdots &
\cos({\nicefrac{1}{2}}P_n) \\
\cos({\nicefrac{3}{2}}P_1) & \cdots & \cos(\nicefrac{3}{2}P_n) \\ 
\vdots & & \vdots \\ \cos((n-\nicefrac{1}{2})P_1) &
\cdots & \cos((n-\nicefrac{1}{2})P_n) \end{pmatrix} .$$
Note that because of \eqref{nyquist2}, we know that $|P_j/2|<\pi/2$ and
hence that ${\cal W}_n$ is regular. We have
$$|{\cal W}_n| =  2^{n(n-1)/2} \prod_{j=1}^n \cos(P_j/2) \prod_{1 \le i <
j \le n} \left( \cos(P_j) - \cos(P_i) \right).$$
Formula \eqref{half} is used to transform the values $\cos(kP_j)$ of 
\eqref{alpha2} into those needed in ${\cal W}_n$. 
Note that for $k=0$, \eqref{half} has $\cos^2(P_j/2)$ in the left-hand
side. Since $|P_j|\le \pi$, we have $|P_j/2|\le \pi/2$ and hence we know
that $\cos(P_j/2) \ge 0$. For $k>0$, \eqref{half} is used by
gradually incrementing $k$ from $1$ to $n-1$.

\section{Numerical illustration}

To illustrate the use of the new matrix pencil formulations for the
trigonometric spectral analysis problems, we consider a texture
denoising and decomposition problem of the form
$$f(x,y) = \sum_{j=1}^n \alpha_j \left( \cos(\phi_j x+ \psi_j y) +
\sin(\phi_j x + \psi_j y) \right),$$
where $f(x,y)$ is only known through its samples, namely the pixel values of the
texture image in Figure 1, 
and not through its function statement, since the parameters $\phi_j, \psi_j, 
\alpha_j, j=1, \ldots, n$ are to be revealed. To solve this problem
we apply the technique explained in Section 5 to the analysis of $f(x,y)$ given
here. The illustration is merely meant as a proof of principle.

\begin{figure}[h!]
\begin{center}
\includegraphics[scale=0.5]{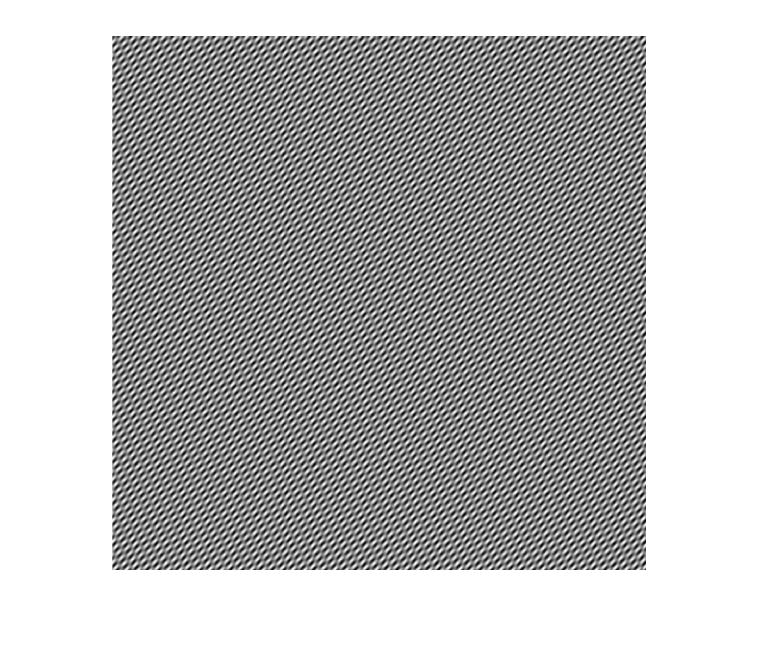}
  \includegraphics[scale=0.5]{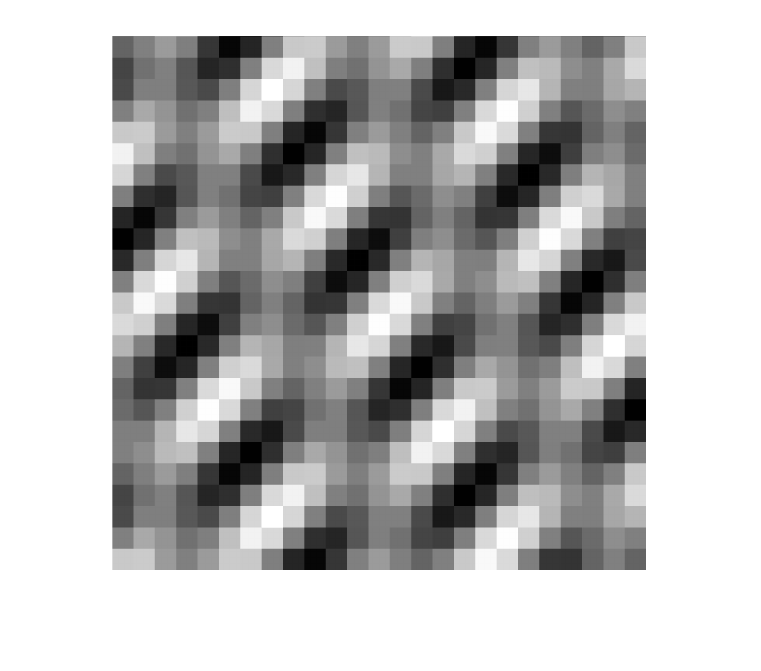}
  \caption{Noise-free texture ($512 \times 512$ left and zoomed-in view on the right). }
  \end{center}
\end{figure}

Note that when using either \eqref{gep_c} or \eqref{gep_s}, the involved 
Hankel and Toeplitz matrices are filled with samples of either the even part
or the odd part of the texture image. So the Hankel and Toeplitz matrices in the
cosine case are not identical to the Hankel and Toeplitz matrices in the
sine case. We proceed with \eqref{gep_c} and hence define 
$$g_k = {1 \over 2} \left( f_k + f_{-k} \right), \qquad 
f_k=f(k\Delta_1, k\Delta_2),$$
in order to obtain the values $\cos(P_j), j=1, \ldots, n$ 
as the generalised eigenvalues 
of the generalised eigenvalue problem
\begin{equation}
\label{texture}
{1 \over 2} \left( H_n^{(1)}+H_n^{(-1)}+T_n^{(1)}+T_n^{(-1)} \right) v_j = 
\cos(P_j)
\left( H_n^{(0)}+T_n^{(0)} \right) v_j. 
\end{equation} 
We further sample
\begin{align*}
F_{\ell,\pm} &= f(\ell\Delta_1\pm\delta_1, \ell\Delta_2\pm\delta_2), \\
G_{k,+} &= {1 \over 2} \left( F_{k,+} + F_{-k,-} \right), \\
G_{k,-} &= {1 \over 2} \left( F_{k,-} + F_{-k,+} \right), 
\end{align*}
to compute the $\cos(Q_j)$.
With $\left({\cal V}_n\right)_{kj} := \cos ((k-1)P_j)$,
we find that the $\cos(Q_j)$, $j=1, \ldots, n$ are obtained from
\eqref{VDM2} and the $\alpha_j$, $j=1, \ldots, n$ from \eqref{alpha2}.
For our image we have $n=2$ and
$$\begin{aligned}
(\phi_1, \psi_1) &= (0.3\pi, 0.15\pi), \\ (\phi_2, \psi_2) &= (0.2\pi,
0.35\pi), \end{aligned}
\qquad \begin{aligned} \alpha_1 &= 1, \\ \alpha_2 &= 1. 
\end{aligned}$$

In practice, the collected data are of course often contaminated with
noise and then the above numerical linear algebra problems are solved in
the least-squares sense. For their rectangular $m \times n$ counterparts
$C_{m,n}^{(r)}$ and $S_{m,n}^{(r)}$ the factorisations given in Theorems 1--4 remain valid, with only the first factor taking size
$m \times n$. Let us add 10 dB of white Gaussian noise to the
image, resulting in Figure 2.

\begin{figure}[h!]
\begin{center}
\includegraphics[scale=0.5]{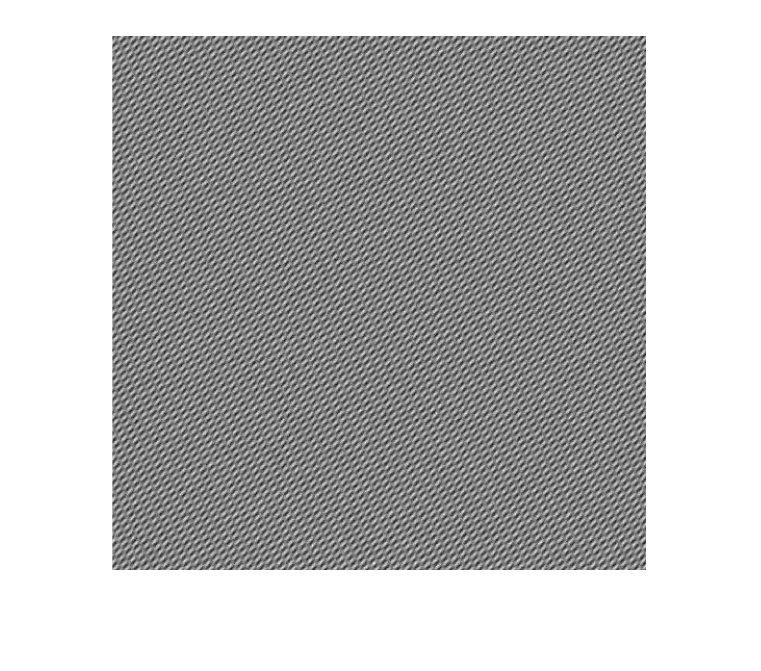}
  \includegraphics[scale=0.5]{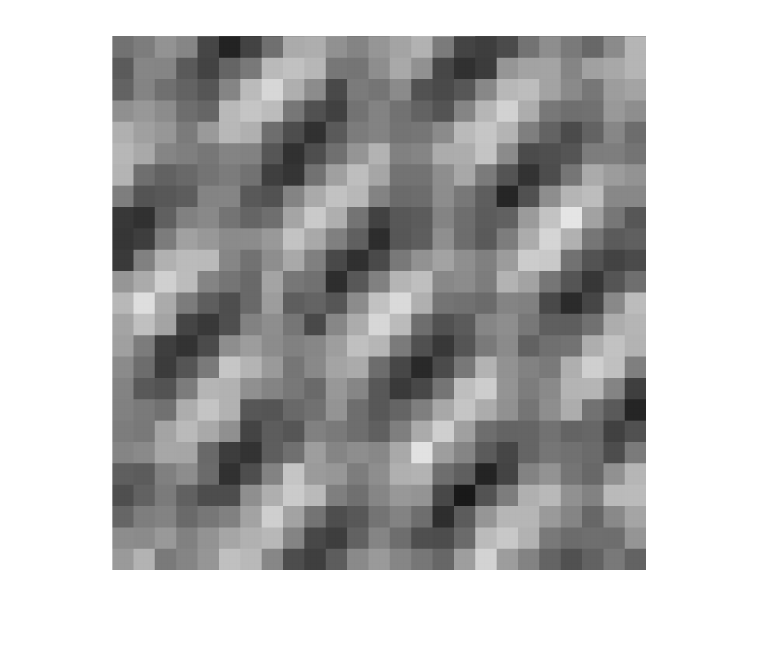}
  \caption{Noisy texture ($512 \times 512$ left and zoomed-in view on the right). }
  \end{center}
\end{figure}

Instead of working with the matrices $H^{(r)}_n$ and $T^{(r)}_n$ in
\eqref{texture} where $r=-1, 0, 1$, we
consider their tall rectangular $m \times \eta$ versions
$H^{(r)}_{m,\eta}$ and $T^{(r)}_{m,\eta}$ having $m$ rows and $\eta$
columns with $m \ge \eta \ge n$. The $\eta-n$ additional terms model the
noise present in the data.

Let the pixels be indexed left to right $-255, \ldots, 0$, $\ldots, 256$
and likewise from top to bottom. We choose $(\Delta_1, \Delta_2) = (1,0)$ and
$(\delta_1, \delta_2)=(0,1)$ and sample the image for simplicity
on 11 full lines (which is much more than we need), namely the
line $y=0$ for the samples $f_k$, and the lines $y=\pm5, \ldots, \pm1$
for the samples $F_{\ell,\pm}$.
These 11 lines constitute 2.15\% of the image. We compose the required
$g_k$ and $G_{k,\pm}$ from the collected pixels. 
We take $m=172$ and $\eta=84$ and reduce the size of each of the Hankel
and Toeplitz matrices involved in the generalised eigenvalue problem
to $m \times 2$ using their singular value decomposition.

The procedure described in Section 5, when applied to the even part of
the image in Figure 1, eventually returns the following values 
\begin{align*}
(\phi_1, \psi_1) &= (0.299993\pi, 0.150008\pi), \\ (\phi_2, \psi_2) &=
(0.200019\pi, 0.349980\pi), \\
\alpha_1 &= 0.9993961, \\ \alpha_2 &= 1.000039 
\end{align*}
which allow us to restore the noisy image; the result, shown in Figure 3, is visually
indistinguishable from the original.
After denoising the following minor discrepancy from the original $f(x,y)$
remains: 
\begin{align*}
&\text{RMSE}_{(\phi,\psi)} = 1.78 \times 10^{-4}, \\
&\text{RMSE}_{\alpha} = 4.28 \times 10^{-4}, \\
&\text{RMSE}_{\text{image}} = 0.0344.
\end{align*}

\begin{figure}[h!]
\begin{center}
\includegraphics[scale=0.5]{image_noisy_zoom}
  \includegraphics[scale=0.5]{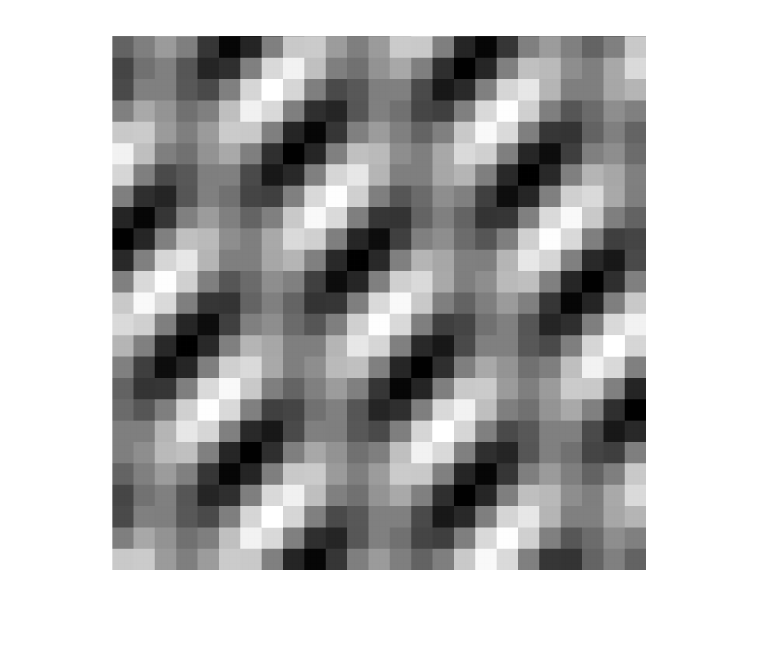}
  \caption{Noisy texture at the left and denoised texture at the right.}
  \end{center}
\end{figure}

\section*{Author contributions}
Both authors contributed equally to this work.

\end{document}